\documentclass[12pt]{amsart}
\usepackage{amscd,amsmath,amssymb,amsfonts}
\renewcommand{\baselinestretch}{1.3}

\makeatletter
\newcommand{\single}{\let\CS=\@currsize\renewcommand{\baselinestretch}{1}\tiny\CS}
\newcommand{\singles}{\let\CS=\@currsize\renewcommand{\baselinestretch}{1.3}\tiny\CS}
\newcommand{\oneanda}{\let\CS=\@currsize\renewcommand{\baselinestretch}{ 1.2}\tiny\CS}
\newcommand{\doubles}{\let\CS=\@currsize\renewcommand{\baselinestretch}{1.5}\tiny\CS}
\newcommand{\tree}{\let\CS=\@currsize\renewcommand{\baselinestretch}{1.5}\tiny\CS}
\newcommand{\four}{\let\CS=\@currsize\renewcommand{\baselinestretch}{2}\tiny\CS}

\newcommand{\ncom}{\newcommand}
\ncom{\bq}{\begin{equation}}
\ncom{\eq}{\end{equation}}
\ncom{\beqn}{\begin{eqnarray*}}
\ncom{\eeqn}{\end{eqnarray*}}
\ncom{\beq}{\begin{eqnarray}}
\ncom{\eeq}{\end{eqnarray}}
\ncom{\been}{\begin{enumerate}}
\ncom{\eeen}{\end{enumerate}}
\ncom{\FP}{{\mathcal H}'_1}
\ncom{\nno}{\nonumber}
\ncom{\ko}{h^0}
\ncom{\Ko}{K^0}
\ncom{\hs}{\mbox{\hspace{.25cm}}}
\ncom{\rar}{\rightarrow}
\ncom{\lrar}{\longrightarrow}
\ncom{\Rar}{\Rightarrow}
\ncom{\noin}{\noindent}
\ncom{\card}{card}
\newtheorem{thm}{Theorem}[section]
\newtheorem{lemma}[thm]{Lemma}
\newtheorem{corollary}[thm]{Corollary}
\newtheorem{pro}[thm]{Proposition}
\newtheorem{example}[thm]{Example}
\newtheorem{remark}[thm]{Remark}
\newtheorem{defn}[thm]{Definition}

\newtheorem{Conj}[thm]{Conjecture}
\ncom{\bc}{\begin{Conj}}
\ncom{\ec}{\end{Conj}}
\ncom{\bd}{\begin{defn}}
\ncom{\ed}{\end{defn}}
\ncom{\bt}{\begin{thm}}
\ncom{\et}{\end{thm}}
\ncom{\dd}{{\bf d}}
\ncom{\bl}{\begin{lemma}}
\ncom{\el}{\end{lemma}}
\ncom{\bsl}{\begin{lemma}}
\ncom{\esl}{\end{lemma}}
\ncom{\bco}{\begin{corollary}}
\ncom{\eco}{\end{corollary}}
\ncom{\bp}{\begin{pro}}
\ncom{\ep}{\end{pro}}
\ncom{\bex}{\begin{example}}
\ncom{\eex}{\end{example}}
\ncom{\brm}{\begin{remark}}
\ncom{\erm}{\end{remark}}
\ncom{\ga}[2]{\begin{gather}\label{\card 1}\card 2 \end{gather}}
\ncom{\llra}{\Longleftrightarrow}
\ncom{\eps}{\epsilon}
\ncom{\comx}{{\mathbb C}}
\ncom{\zee}{$Z\!\!\!\!Z$}
\ncom{\ze}{{\mathbb Z}}
\ncom{\E}{{\mathbb E}}
\ncom{\Q}{{\mathbb Q}}
\ncom{\p}{{\mathbb P}}
\ncom{\pn}{{\mathbb P}^n}
\ncom{\R}{{\mathbb R}}
\ncom{\lla}{{\lambda}}
\ncom{\lls}{{\ell^*}}
\ncom{\llp}{{\ell'^*}}
\ncom{\LL}{{\mathcal L}}
\ncom{\hLL}{\hat{\mathcal L}}
\ncom{\TT}{{\mathcal T}}
\ncom{\UU}{{\mathcal Y}}
\ncom{\hUU}{\hat{\mathcal Y}}
\ncom{\HHH}{{\mathcal H}}
\ncom{\HHHH}{H}
\ncom{\KKKK}{L'}
\ncom{\LLLL}{L}
\ncom{\PP}{P}
\ncom{\PPP}{{P}}
\ncom{\DDD}{{\mathbb D}}
\ncom{\DD}{{\mathcal D}}
\ncom{\un}{{\bf 1}}

\ncom{\deu}{{\bf 2}}

\ncom{\ZZ}{Z}
\ncom{\OO}{{\mathcal O}}
\ncom{\FF}{{\mathcal F}}
\ncom{\KK}{{\mathcal K}}
\ncom{\KKK}{{\mathcal K}}
\ncom{\HH}{{\mathcal H}}
\ncom{\QCH}{{\rm QCH}}
\ncom{\FH}{FH}
\ncom{\K}{{\mathbb K}}
\ncom{\AAA}{{\mathbb A}}
\ncom{\G}{{\mathbb G}}
\ncom{\D}{{\mathbb D}}
\ncom{\Y}{{\mathbb Y}}
\ncom{\F}{{\mathbb F}}
\ncom{\h}{{\mathbb H}}
\ncom{\al}{\alpha}
\ncom{\be}{\beta}
\ncom{\f}{\frac}
\ncom{\gam}{\gamma}
\ncom{\bib}{\bibitem}
\ncom{\pf}{{\bf Proof: }}
\ncom{\sta}{\stackrel}
\ncom{\ov}{\overline}
\ncom{\cX}{{\mathcal X}}
\ncom{\cA}{{\mathcal A}}
\ncom{\cB}{{\mathcal B}}
\ncom{\cG}{{\mathcal G}}
\ncom{\cI}{{\mathcal I}}
\ncom{\cO}{{\mathcal O}}
\ncom{\cV}{W}
\ncom{\cW}{{\mathcal W}}
\ncom{\cK}{{\mathcal K}}
\ncom{\cE}{{\mathcal E}}
\ncom{\cL}{{\mathcal L}}
\ncom{\cZ}{{\mathcal Z}}
\ncom{\cP}{{\mathcal P}}
\ncom{\cN}{{\mathcal N}}
\ncom{\cF}{{\mathcal F}}
\ncom{\cM}{{\mathcal M}}
\ncom{\cC}{{\mathcal C}}
\ncom{\cR}{{\mathcal R}}
\ncom{\cS}{{\mathcal S}}
\ncom{\cH}{{\mathcal H}}
\ncom{\cY}{{\mathcal Y}}
\ncom{\s}{(\!\!\times}
\ncom{\m}{\mbox}
\ncom{\eop}{{\hfill $\Box$}}

\ncom{\Ho}{H^0}
\ncom{\ho}{h^0}
\ncom{\ttt}{\mathbf{t}}
\ncom{\ee}{\mathbf{e}}

\title[Hilbert schemes of fat $r$-planes]{Hilbert schemes of fat $r$-planes and the triviality of\\ Chow groups
of complete intersections }
\author[A.Hirschowitz]{Andr\'e Hirschowitz }
\address{CNRS, Laboratoire J.-A.Dieudonn\'e, Universit\'e de Nice--Sophia Antipolis,
Parc Valrose, 06108 Nice Cedex 02, France}
\email{ah@math.unice.fr}

\author[J. N. Iyer]{Jaya NN Iyer}

\address{The Institute of Mathematical Sciences, CIT
Campus, Taramani, Chennai 600113, India}
\email{jniyer@imsc.res.in }

\begin{document}

\footnotetext{Mathematics Classification Number: 14C25, 14D07, 14D22, 14F40.}
\footnotetext{Keywords: Complete intersections, Chow groups, Hilbert schemes.}
\begin{abstract}
In this paper, we investigate the question of triviality of the rational Chow groups of complete intersections
in projective spaces and obtain improved bounds for this triviality to hold. Along the way, we have to
study the dimension and nonemptiness of some Hilbert schemes of fat $r$-planes
contained in a complete intersection $Y$, generalizing well-known 
results on the Fano varieties of $r$-planes contained in $Y$.
\end{abstract}

\maketitle

\setcounter{tocdepth}{1}
\tableofcontents


\section{Introduction}


\subsection{The triviality conjecture}

The aim of this paper is to investigate the triviality of the
low-dimensional rational Chow groups for certain projective varieties.
If $Y$ is a nonsingular complete intersection of multidegree $(d_1, \cdots ,d_s)$
in a projective space $\p^n$, and $n$ is sufficiently large with respect to the degrees,
it is known that, for small values of $r$, the rational Chow group $QCH_r(Y):= CH_r(Y)\otimes \Q$ is trivial,
namely one-dimensional (generated by the linear sections).
The precise conjectural bound on the multidegrees for the triviality follows from the
study of the \textit{Hodge type} of the complementary open variety $\p^n-Y$ initiated
by Deligne \cite{De} and followed by works of Deligne-Dimca \cite{De-Di} and
Esnault-Nori-Srinivas \cite{Es}, \cite{ENS}. A formulation of the conjectured bound was
made in \cite[Conjecture 1.9]{Pa}, which says:

\bc\label{conjecture}
Suppose $Y\subset \p^n$ is a smooth complete intersection of
multidegree $(d_1, \cdots ,d_s)$, and let $r$ be a nonnegative integer.
If
$$
r(\m{max}_{1\leq i\leq s}d_i)+\sum_{i=1}^sd_i\leq n
$$
then $\QCH_r(Y)$ is trivial.
\ec

This conjecture entails the hypersurface case ($s=1$) as well as the higher-dimensional case $s \ge 2$, and, as we will see, our 
contribution concerns mostly the latter.

\subsection{Fat and strong planes}

The central notions in our approach are those of {\em fat} and {\em strong} planes, which appear at least 
implicitly in \cite{ELV}, and go back to Roitman for the $0$-dimensional case.

By a $t$-fat $r$-plane in a projective space, we mean the $t$-th infinitesimal neighborhood of an
$r$-plane in an $(r+1)$-plane. Given a subscheme $Y'$ in a projective space and a Cartier divisor $Y$
in $Y'$, we say that an $r$-plane $L$ in $Y$ is strong (with respect to $Y'$) if there exists an $(r+1)$-plane
$L'$ in $Y'$ containing $L$ such that the set-theoretic intersection $L' \cap Y$ is either $L$ or $L'$.
The connection between the two notions is given by the following statement, proven in section \ref{strongplanes} :

\bp 
Suppose $Y$ is in the linear system $| {\OO}_{Y'}(t)|$. Then any strong $r$-plane in $Y$ is the support of a
$t$-fat $r$-plane contained in $Y$. Conversely, if furthermore $Y'$ is (set-theoretically) defined by
equations of degree (strictly) less than $t$, then the support of any $t$-fat $r$-plane contained in $Y$ is strong.
\ep

\subsection{Roitman's technique, small steps and big steps}

The case of $0$-cycles has been handled by Roitman \cite{Ro}. His method consists in starting from a
(positive) $0$-cycle $Z$ on $Y$, and building a \textbf{ruled} cycle in $\p^n$ whose intersection with $Y$ will be 
not too far from a multiple
of $Z$. This is achieved by choosing a ruling by lines which are \textbf{strong} (see below), in the sense 
that either they cut $Y$ in a single (multiple) point or they are inside $Y$. This method can be extended to the higher-dimensional case of 
$r$-cycles. Of course, the scope of this method is limited by the need of ``sufficiently many strong $(r+1)$-planes''.
This approach has been successfully applied in\cite{ELV} through a single big step, showing that the 
restriction 
of $(r+s)$-cycles from $\p^n$ to $Y$ is sufficiently surjective under the following numerical assumption (at least for 
degrees at least $3$, the assumption being different when all degrees are equal to $2$):
$$
n \ge \sum_{i=1}^s{d_i+r \choose r+1}.
$$

The geometric meaning of their numerical condition is the rational-connectedness
of the variety of $r$-planes in $Y$.

In the present work, we apply Roitman's technique through smaller steps, typically showing that the restriction of $(r+1)$-cycles on the suitable complete intersection of
multidegree $(d_1, \cdots ,d_{s-1})$ to $Y$ is sufficiently surjective. The analysis of small steps being somewhat simpler, we succeed in 
applying Roitman's technique to small steps in essentially the whole expected range, relaxing in particular the rational-connectedness assumption.
For instance for $5$-cycles on complete intersections of type $(20, 30)$, the rational-connectedness condition requires $n \ge 1 800 000$ while we take care of all cases with  $n \ge 370 000$.

\subsection{Our small step theorem}

Our first main result reads as follows:

\bt
\label{Surj}
If, in the subvariety $ Y' \subset \p^n$, the Cartier divisor $Y \in |\OO_{Y'}(t)|$ is covered by strong $r$-planes,
then
the restriction map
$ {\rm QCH}_{r+1} (Y') \rar {\rm QCH}_r (Y) $ is onto.
\et

This implies in particular that whenever $ {\rm QCH}_{r+1} (Y')$
is trivial, so is $ {\rm QCH}_r (Y) $.

Our proof of Theorem \ref{Surj} follows the corresponding proof for the big step in \cite{ELV}, with a single but decisive
technical improvement: we introduce a different filtration of the Chow group $CH_r$, where $CH_r^{(s)}$ is generated by subvarieties 
covered by \textbf{strong} $s$-planes. We did not explore yet whether such a new filtration could also improve the bound for the big step.

In order to apply the above result, we need to find the appropriate condition on
the degrees for our
complete intersection $Y$ to be covered by strong $r$-planes.

\subsection{Covering by strong planes}

Thus we are led to search for the numerical condition for at least the generic 
complete intersection of multi-degree $(d_1, \cdots , d_s)$ to be covered by strong $r$-planes.

We recall that the strongness property is with respect to the pair $(Y, Y')$.
We say that an $s$-codimensional subvariety in a projective space has type $(d_1, \cdots , d_s)$ when it
is a union of irreducible components of a complete intersection of multi-degree $(d_1, \cdots , d_s)$.
Accordingly, we say that a pair
$Y \subset Y'$ has type $(d_1, \cdots , d_s)$ if $Y'$ has type $(d_1, \cdots , d_{s-1})$ and $Y$ is a divisor of
degree $d_s$ in $Y'$. In \S \ref{spannedness}, we prove :

\bp
\label{Covs}
Let $n, r, s, d_1, \cdots , d_s$ be integers satisfying
$$r \ge 0, \,~~ s \ge 1, ~~\, n \ge r+s, \,~~ 2 \le d_1 \le \cdots \le d_{s-1} < d_s$$
and the ("expected") inequality

\vspace {-0.6cm}
$$\rho+r \ge n-s$$

where

\vspace {-0.6cm}

$$\rho := (r+2)(n-r) -\sum_{i=1}^s {d_i + r +1 \choose r+1}$$

is the  dimension of the variety of $d_s$-fat $r$-planes in the general
complete intersection of type $(d_1, \cdots , d_s)$ (see \S \ref{fatplanes}).

If $Y \subset Y'$ is any pair of type $(d_1, \cdots , d_s)$ in $\p^n$, then $Y$ is covered
by strong $r$-planes.
\ep

Note that the intended meaning of the "expected" inequality relates the dimension of the
universal
$d_s$-fat $r$-plane with the dimension of our complete intersection $Y$.

Also note the strict inequality
$d_{s-1} < d_s$. Apart from this restriction, our result is the expected one. The
discarded case
would involve a refined analysis (this is where we do not cover the whole
range of Roitman's method for small steps).

Our proof of Proposition \ref{Covs} relies on the study of the Hilbert schemes of fat $r$-planes contained in a general 
complete
intersection.
We show in \S \ref{maxrank} that they have the expected dimension; but we need a more accurate result saying
that, when this expected dimension is
nonnegative, these Hilbert schemes are nonempty. We conjecture that this is true in most cases (despite the 
notable
exception of double lines on quadric surfaces), and prove it in the case we need for our application to Chow
groups. For such a result, as illustrated in \cite{De-Ma}, two approaches are available:
through intersection computations or through maximal rank problems. We follow the latter approach, using a method 
that can be tracked back  at least to \cite{EH, EHM}. 

\subsection{The main theorem}

Combining the previous
results,
we obtain our main result:

\bt
\label{Coro}
Let $n, r, s, d_1 \le \cdots \le d_{s-1} < d_s$ be integers as above, satisfying 
$$\rho+r \ge n-s.$$

If $Y \subset Y'$ is any pair of type $(d_1, \cdots , d_s)$ in $\p^n$, and if $ {\rm QCH}_{r+1} (Y')$
is trivial, then so is $ {\rm QCH}_r (Y) $.
\et

This theorem may be applied recursively. For instance in the case of codimension two complete intersections ($s=2$),
our assumption for triviality reads, for ($d_1 <d_2$) :

$$ n \ge \text{max} ( \frac{{d_1+r+1 \choose r+1}+  {d_2+r+1 \choose r+1}+ r^2+r-2 }{r+1}, 
 \frac{{d_1+r+2 \choose r+2}+ r^2+3r}{r+2}  ).$$

In order to compare this new bound with \cite {ELV}'s

$$
n \ge  {d_1+r \choose r+1}+  {d_2+r \choose r+1}
$$

 we fix $r$ and $d_1$ and let  $d_2$ vary. For sufficiently large values of $d_2$, the new bound can
 be estimated as $\frac {d_2^{r+1}}{(r+1) (r+1)!}$. While if one uses \cite{ELV}, the best value of $n$ is
estimated as
$\frac {d_2^{r+1}}{(r+1)!}$. Hence in this context of a large $d_2$, we roughly divide by $r+1$ the range
where the conjecture is still
open. On the other hand, it can be checked than if $d_2$ is not sufficiently large, \cite {ELV}'s bound
 remains the best (for instance fix $d_2 =d_1+1$ and let $d_1$ go to infinity).

In the higher-codimensional cases, a similar picture will occur, namely our result will 
provide an improved bound only for sufficiently large values of $d_s$. Furthermore, in many cases, the best
bound will be obtained by combining one or more of our small steps with a big step from \cite{ELV}.

\subsection{The case of hypersurfaces}

The case of hypersurfaces ($s=1$) has been considered in the first place.
Concerning a general cubic hypersurface $Y\subset \p^n$, C. Schoen \cite{Sc} showed the
triviality $\QCH_1(Y)\simeq \Q$ when $n\geq 7$ and Paranjape \cite{Pa} obtained the sharp
bound in this case showing the triviality of $1$-cycles when $n\geq 6$ (in the same paper, he gave
the first finite bound
for general complete intersections). 

For hypersurfaces of the general degree $d$, the best known bound has been obtained "in the margin" by 
J. Lewis \cite{Le2} (added
on proofs at the very end of the paper). There, the statement concerns only the generic hypersurface, and
the bound occurs as the condition for the so-called
cylinder homomorphism to be surjective. This bound by Lewis is better than the bound obtained later 
(for the hypersurface case)
in \cite{ELV}. It was 
rediscovered by A. Otwinowska \cite{Ot}: there the statement concerns all smooth hypersurfaces, and the 
geometric meaning of the bound is that the  hypersurfaces
of degree $d$ in $\p^{n+1}$ (not $\p^n$!) are covered by $(r+1)$-planes. 
Surprisingly, our small step gives exactly the same bound, with a third geometric meaning for the condition, 
namely that the 
hypersurfaces are covered by $d$-fat $r$-planes. Furthermore, our statement concerns all  hypersurfaces, not only 
smooth ones.

\subsection{The base field}

We work over an algebraically
closed field of characteristic zero.
The closedness assumption could be removed, thanks to the fact that the kernel of 
${\rm CH}_r(Y_k) \to {\rm CH}_r(Y_{\bar{k}})$ is
torsion \cite {Bl2}, while
the characteristic zero assumption is used in the proofs of \S \ref {fatplanes}.

{\Small
\textit{Acknowledgements}: It is a pleasure to dedicate this paper in honour of S. Ramanan. Both authors have
experienced very fruitful mathematical interaction with him, and take the opportunity to
acknowledge his deep influence.
 This work was initiated during the second author's stay at MPI, Bonn in 2003 and 
partly done during her visit to Nice in Dec 2004 and at IAS, Princeton in 2007. The support and hospitality of 
these institutions is gratefully acknowledged. We also thank H. Esnault, J. Lewis and M. Levine for their feedback at some point during the course of this work.
}


\section{Strong planes}\label{strongplanes}


Throughout this section, we consider a subvariety $Y'$ in a projective space
equipped with a Cartier divisor $Y$, and we fix an integer $r$. We are interested in the restriction map
$ {\rm QCH}_{r+1} (Y') \rar {\rm QCH}_r (Y) $, where we write $ {\rm QCH}_{r} (W) $
for the rational Chow group of $r$-dimensional cycles on $W$.

Recall that an $r$-plane $L$ in $Y$ is said strong (with respect to $Y'$) if there exists an $(r+1)$-plane
$L'$ in $Y'$ containing $L$ such that the set-theoretic intersection $L' \cap Y$ is either $L$ or $L'$.

In this section, we prove our first main result :

\bt
If, in the subvariety $ Y' \subset \p^n$, the Cartier divisor $Y \in |\OO_{Y'}(d)|$ is covered by strong $r$-planes,
then
the restriction map
$ {\rm QCH}_{r+1} (Y') \rar {\rm QCH}_r (Y) $ is onto.
\et

For the proof, we generalize our notion of strongness and define a notion of strong $s$-plane in $Y'$ for
$s \le r+1$. A $(r+1)$-plane $H$ in $Y'$ is said {\it strong} (with respect to the pair $(Y, Y')$) if it is 
contained in $Y$,
or if its set-theoretic intersection
with $Y$ is a $r$-plane. Then, for $s \le r$, a $s$-plane in $Y'$
is said to be {\it $r$-strong}, or simply {\it strong} (when $r$ is clear from the context), if it is contained in a strong $(r+1)$-plane.
As usual, we say that a closed subvariety $W$ of $Y$ is
{\it spanned} or {\it covered} by strong $s$-planes if it is a union of strong $s$-planes contained in $W$.

Now we denote by ${\rm QCH}^{(s)}_{r} (Y)$ the subgroup of ${\rm QCH}_r (Y)$
which is generated by $r$-dimensional subvarieties of $Y$ which
are spanned by strong $s$-planes. This is the place where our proof differs from the corresponding proof in \cite {ELV}.
Note that any subvariety in $Y$ is spanned at least by strong
$0$-planes: since $Y$ is covered by strong $r$-planes, it is also covered by strong $0$-planes.
Thus we have ${\rm QCH}^{(0)}_{r} (Y) = {\rm QCH}_r (Y)$.

For $s \ge 1$, if $Z$ is spanned by strong $s$-planes it is spanned by strong $(s-1)$-planes
as well. Hence one has ${\rm QCH}^{(s)}_{r} (Y) \subseteq
{\rm QCH}^{(s-1)}_{r} (Y)$. For $s > r$ one has ${\rm QCH}^{(s)}_{r} (Y) = \{ 0 \}$.

We prove by descending induction on $s$ that ${\rm QCH}^{(s)}_{r} (Y)$ is in the image of
${\rm QCH}_{r+1} (Y')$. The initial case is with $s :=r+1$ and follows since
 ${\rm QCH}^{(r+1)}_{r} (Y)$ is reduced to $0$.

Before stating the induction step as a lemma, we introduce the following notation.
Let $\Gamma \subset Y'$ be an $(r+1)$-dimensional closed
subvariety or, more generally, an $(r+1)$-cycle. By \cite[8.1]{Fu}, the intersection product $\Gamma
\cdot Y$ is a class in ${\rm CH}_r (|\Gamma | \cap Y)$. By abuse of
notation we will also write $\Gamma \cdot Y$ for its image in  ${\rm QCH}_r(Y)$.

\bl
\label{Pro1}
Let $s$ be an integer with $0 \le s \le r$, and $W$ be an $r$-dimensional irreducible subvariety of $Y$, 
spanned
by strong $s$-planes but not by strong $(s+1)$-planes.
Then there exist an $(r+1)$-dimensional cycle $\Gamma$ in
$Y'$ and a positive integer $\alpha$ with
$$
\Gamma \cdot Y \equiv \alpha W \ {\rm mod} \ {\rm QCH}^{(s+1)}_{r} (Y).
$$
\el
\begin{proof}
We start with the case $s:=r$ which means that $W$ is a strong $r$-plane. This gives us a strong
$(r+1)$-plane in $Y'$ which we take for $\Gamma$. Indeed, we have 
$\Gamma \cdot Y = d W .$

Now we suppose $s < r$.
In order to define $\Gamma$, we start by choosing carefully an algebraic family
$(H_z)_{z\in Z}$ of strong $s$-planes covering $W$. Note that by our assumption on $s$,
each strong $s$-plane in $Y$ is contained in a strong $(s+1)$-plane also contained in $Y$, thus we may choose 
more precisely
an algebraic family $(H_z \subset H'_z)_{z\in Z}$ where $(H'_z)$ is a strong $(s+1)$-plane in $Y$,
$H_z$ is a hyperplane in $H'_z$ and $W$ is covered by $(H_z)_{z\in Z}$. By standard arguments, we may suppose 
that $Z$ is projective smooth connected of dimension $r-s$. We denote by $H_Z \subset H'_Z$ the two corresponding projective bundles 
over $Z$.

Since $W$ is not covered by strong $(s+1)$-planes, the projection of $H'_Z$ into $Y$ is not contained in $W$, 
thus it is a positive $(r+1)$-cycle. We take for $\Gamma$ this Chow-theoretic projection of $H'_Z$ in $Y'$.

Let us now compute $\Gamma \cdot Y$ in ${\rm QCH}_{r} (Y)$ (remind that we consider $\Gamma$ as a cycle in $Y'$).

We start by applying the projection formula
\cite[8.1.7]{Fu} to $pr_2:Z \times Y' \to
Y'$:

$$
\Gamma \cdot Y = pr_{2*} (H'_Z) \cdot Y = pr_{2*}
(H'_Z \cdot (Z \times Y)).$$

So now we compute $H'_Z \cdot (Z \times Y).$ This is the divisor
class in $H'_Z$ defined by
the linear system $|pr_2^*\OO_{Y'}(Y)|$.
Now $H'_Z$ is a projective bundle and this linear system has degree $d$ along the fibers of this bundle. 
Thus we have

$$ H'_Z \cdot (Z \times Y) = d H_Z + \psi^{-1} (D)
$$

where  $D$ is a divisor in $Z$ and $\psi : H'_Z \to Z$ is the bundle projection. 
We get

$$ \Gamma \cdot Y = d pr_{2*} (H_Z) + pr_{2*} \psi^{-1} (D)
$$
in ${\rm QCH}_r (Y \cap pr_2 (H'_Z ))$. Since $H_Z$ is
generically finite over the subvariety $W$ and since $pr_{2*}
(\psi^{-1} (D))$ lies in ${\rm QCH}^{(s+1)}_{r} (Y)$, one
obtains, for some positive multiple $\alpha$ of $d$, the relation
$$
\Gamma \cdot Y \equiv \alpha W \ {\rm mod} \ {\rm QCH}^{(s+1)}_{r} (Y).
$$
\end{proof}

Now we check the following statement,
already mentioned in our introduction:

\bp \label{fatstrong}
Suppose $Y$ is in the linear system $| {\OO}_{Y'}(t)|$. Then any strong $r$-plane in $Y$ is the support of a
$t$-fat $r$-plane contained in $Y$. Conversely, if furthermore $Y'$ is (set-theoretically) defined by
equations of degree (strictly) less than $t$, then the support of any $t$-fat $r$-plane contained in $Y$ is strong.
\ep
\begin{proof} 
For the first statement, our strong $r$-plane $L$ is contained in a strong $(r+1)$-plane $L'\subset Y'$. If $L'$ 
is contained 
in $Y$, then so is the $t$-th infinitesimal neighborhood of $L$ in $L'$. If not, then, since the set-theoretic 
intersection
of $Y$ and $L'$ is $L$, and the degree of the restriction of $| {\OO}_{Y'}(t)|$ to $L'$ is $t$, the 
scheme-theoretic intersection
$Y \cap L'$ has to be the $t$-th infinitesimal neighborhood of $L$ in $L'$. Thus in both cases, this is 
the desired $t$-fat $r$-plane.

For the second statement, let $L \subset \pn$ be a $t$-fat $r$-plane
contained in the $(r+1)$-plane $L'$.
The equations defining $Y'$ vanish on $L$. Since these equations can be chosen
of degree strictly less than $t$,
they vanish identically on $L'$, which means that $L'$ is contained in $Y'$,
hence that $L$ is strong.

\end{proof}


\section{Restricted flag-Hilbert schemes}\label{hilb}


In this section, we collect some technical material concerning the infinitesimal theory of
restricted flag-Hilbert schemes.
Here by a full Hilbert scheme (for a given projective variety), we mean any open subscheme of the  Hilbert scheme
associated to a Hilbert polynomial, while by a Hilbert scheme, we mean any locally closed subscheme of a full Hilbert scheme.

Given two full Hilbert schemes $\HHH _1$ and $\HHH_2$ of subschemes of the same ambient
projective scheme $P$, we have the corresponding flag-Hilbert scheme $\DD$ of pairs $(X \hookrightarrow Y)$ in
$\HHH_1 \times \HHH_2$.
Two subschemes $ \HHH' _1 \subset \HHH _1 $ and $\HHH' _2 \subset \HHH_2$ being given,
by the corresponding \emph {restricted} flag-Hilbert scheme,
we mean the scheme-theoretic intersection $\DD'$ of $\DD$ with $\HHH'_1 \times \HHH'_2$.

In the example we have in mind, $P$ is a projective space, 
$\HHH'_1$ is a variety of fat planes, and $\HHH'_2 = \HHH_2$ is a full Hilbert scheme
of complete intersections.

We write $i : X \to Y$ for a given pair, $ \mathcal{ I} _X$ and $ \mathcal{ I} _Y$ for the two ideal sheaves on $\PPP$,
$N_X:= \mathcal{ H}om ( \mathcal{ I} _X, \OO_X) $ and $N_Y:= \mathcal{ H}om ( \mathcal{ I} _Y, \OO_Y) $
for the corresponding normal bundles.
We denote by $ N_{Y|X}$ the restriction $i^*(N_Y)$ of $N_Y$ to $X$.
We also have $ i_*: N_X \to  i^*N_{Y} $ and
 $i^* : N_Y \to i_*  N_{Y|X} $. Note that the two codomains have the same space of sections $H^0(N_{Y|X} )$.
Putting together, we have a morphism 
$$
(i_*, i^*) : H^0(N_X) \oplus H^0 (N_Y) \to H^0(N_{Y|X} ).
$$
The domain of this morphism is the tangent space to the product of our two Hilbert schemes, and the tangent space
to the flag-Hilbert scheme is identified as the kernel of the above map $(i_*,i^*)$ (see \cite {Kleppe},
\cite[Remark 4.5.4 ii]{Se}). Hence the differentials of the two
projections are the restrictions  to this kernel of the projections.

We first state in our way the standard result in the unrestricted case:


\bp
We suppose that $\HHH_1, \, \HHH_2$ are smooth connected and that $\DD$ has codimension $c$ at $O :=(X, Y)$.
We also suppose that $ i^* :  H^0 (N_Y) \to H^0(N_{Y|X} )$ has rank $c$. Then

(i) $\DD$ is smooth at $O$;

(ii) The image of $i_* : H^0(N_X) \to H^0(N_{Y|X} )$ is contained in the image of
$i^*:  H^0 (N_Y) \to H^0(N_{Y|X}) $;

(iii) The first projection $\DD \to \HHH_1$ is smooth  at $(X, Y)$;

(iv) The second projection $\DD \to \HHH_2$ is smooth  at $(X, Y)$ if (and only if) the rank of
$i_* : H^0(N_X) \to H^0(N_{Y|X} )$ is $c$.

\ep
\begin{proof}
(i) Since $i^*$ has rank $c$, the pair $(i_*, i^*)$ has rank at least $c$. It follows that, in the tangent space
of $\HHH_1 \times \HHH_2$ at $O$, the tangent space to $\DD$ is at least $c$-codimensional. Since $\DD$ is
$c$-codimensional, this implies that $\DD$ is smooth at $O$.

(ii) By the previous argument, we see that the rank of the pair $(i_*, i^*)$ is exactly $c$, which means
the stated inclusion.

(iii) Using the previous item and an easy diagram-chasing, we see
 that the differential of $\DD \to \HHH_1$ at $O$ is onto.

(iv) This follows by a similar diagram chasing.

\end{proof}

Now we turn to the restricted case. Here we write $\Ko (N_X) \subset H^0(N_X)$ for the tangent space to
$\HHH'_1$ at $X$, and $\Ko (N_Y) \subset H^0(N_Y)$ for the tangent space to
$\HHH'_2$ at $Y$.

\bp \label{restr}
We suppose that $\HHH'_1, \, \HHH'_2$ are smooth connected and that $\DD'$ has codimension $c$ at $O :=(X, Y)$.
We also suppose that $ i^* :  \Ko (N_Y) \to H^0(N_{Y|X} )$ has rank $c$. Then

(i) $\DD'$ is smooth at $O$;

(ii) The image of $i_* : \Ko(N_X) \to H^0(N_{Y|X} )$ is contained in the image of
$i^*:  \Ko (N_Y) \to H^0(N_{Y|X}) $;

(iii) The first projection $\DD' \to \HHH'_1$ is smooth  at $(X, Y)$;

(iv) The second projection $\DD' \to \HHH'_2$ is smooth  at $(X, Y)$ if (and only if) the rank of
$i_* : \Ko(N_X) \to H^0(N_{Y|X} )$ is $c$.

\ep
\begin{proof}
The main point is the identification of the tangent space to $\DD'$: a pair $(t_1, t_2)$
of vectors in
$\Ko(N_X) \times \Ko (N_Y)$ is tangent to $\DD'$ if the subscheme (over Spec $k[\epsilon]$)
corresponding to $t_1$ is included in the one corresponding to $t_2$. This means exactly that $(t_1, t_2)$ 
is tangent to $\DD$. Hence
the tangent space to $\DD'$ is the kernel of the restriction 
$(i_*, i^*) : \Ko(N_X) \oplus \Ko (N_Y) \to H^0(N_{Y|X} )$. The rest of the proof is identical to the previous one.

\end{proof}


\section{Fat planes in complete intersections}\label{fatplanes}

In this section, we consider 
\begin {itemize}
 \item 
a projective space $\pn$, 
\item
an integer $r$ with
$0 \le r < n$, which is the dimension of our (fat) planes,
\item
an integer $s$ with $1 \le s \le n -r -1$, which is the codimension of our complete intersections (or the number of their equations),
\item
a sequence $\dd := (d_1, \cdots, d_s)$ of $s$ positive integers, which is the multidegree of our complete intersections,
\item 
an integer $t$, with $2 \le t \le \max \dd $, which is the multiplicity of our fat $r$-planes.
\end {itemize}

We keep the notations of the previous section for our case where $\HHH'_1$ is the (smooth) Hilbert scheme parametrizing
$t$-fat $r$-planes in $\pn$ and $\HHH'_2 =\HHH_2$
is the (smooth) Hilbert scheme of complete intersections of
type $\dd$. We write $\delta'_i$  for the dimension of $\HHH'_i$.
The dimension $\delta'_1$ of $\HHH'_1$ does not depend on $t$ (thanks to the assumption $t \ge 2$), it is the
dimension of the corresponding flag variety, namely $(r+2)(n-r-1) + r+1$, in other words $(r+2)(n-r)-1$.

We set
$\rho :=
(r+2)(n-r)-1 -\Sigma_{i=1}^s {d_i + r +1 \choose r+1}+\Sigma_{d_i \ge t} {d_i -t + r +1 \choose r+1}$.
We will see that $\rho$ is the expected dimension for the Hilbert scheme of $t$-fat $r$-planes in a complete
intersection of type $\dd$ in $\p^n$. Recall that by a $t$-fat $r$-plane, we mean the $t$-th infinitesimal
neighborhood of an $r$-plane in an
$(r+1)$-plane. Finally, we set $c := \delta'_1 - \rho$. Hence we have $\rho =\delta'_1 -c$ which means that $c$
is the (expected) number of conditions imposed to a $t$-fat $r$-plane for being contained in a given complete
intersection of type $\dd$. The first result of this section confirms this expectation.

\bp
\label{dimfat}

(i) The codimension of the restricted flag-Hilbert scheme $\DD'$ in $\HHH'_1 \times  \HHH'_2$ is $c$;

(ii) For the generic complete intersection $Y$ of type $\dd$ in $\p^n$ the dimension of
the Hilbert scheme of
$t$-fat $r$-planes in $Y$ is everywhere $\rho$. In particular this Hilbert scheme is empty if
$\rho$ is negative.
\ep
\begin{proof}
(i)
We consider a variety $V$ parameterizing our complete intersections, namely the open subset of the vector space
$\hat V$
of $s$-tuples of
homogeneous polynomials (in $n+1$ variables) of the given multidegree defining an $s$-codimensional subscheme in
$ \p^n$. We write $d_V$ for the dimension of $V$. 
 The variety $V$ comes equipped with the tautological subscheme $T \subset V \times \p^n$. The
corresponding morphism $V \to \HHH_2$ is surjective and it is easily checked to be smooth.
Similarly, $\HHH'_1$
comes equipped (thanks to $t \ge 2$) with a tautological flag $\LL \subset \LL' \subset \HHH'_1 \times \p^n$,
where $\LL$ is
the tautological $t$-fat $r$-plane, while $\LL'$ is its linear span: its fibers over $\HHH'_1$ are $(r+1)$-planes.
Next, we introduce the incidence subscheme $D :=  \DD' \times _{\HHH_2}V \subset \HHH'_1  \times V$. 
Since  $V \to \HHH_2$ is surjective and smooth, it is enough to prove that the codimension of $D$ in
$\HHH'_1 \times V$ is $c$.

  The dimension $d_D$
of $D$ is
understood through the
projection on $\HHH'_1$. Indeed the fibers of the projection $D \rar \HHH'_1$ are traces on $V$ of sub-vectorspaces in
$\hat V$. So we have to compute the codimension in $\hat V$ of tuples vanishing on a fixed $t$-fat $r$-plane
$L$. This
codimension is $\Sigma_{i=1}^s c_i$, where $c_i$ is the codimension of homogeneous polynomials of degree $d_i$
vanishing on $L$. We easily check $c_i = {d_i + r +1 \choose r+1} - {d_i -t + r +1 \choose r+1}$, where we adopt the  convention that
${p \choose q}$ is zero whenever $p < q$.
Hence we end up with the desired result
$$
d_D = d_V +\rho.
$$

(ii) This is an immediate consequence of the first item.

\end{proof}

We need a complementary statement which is a particular case of the following conjecture:

\bc
\label{cdimfat}
Apart from the exception below, for the generic complete intersection $Y$ of type $\dd$ in $\p^n$, when $\rho$
is nonnegative,
the Hilbert scheme of
$t$-fat $r$-planes in $Y$ is nonempty.
\ec

Here is the known exception :

\bex
For double lines on the generic quadric in $\p^3$, we have $\rho =0$ while
the corresponding Hilbert scheme  is empty.
\eex

In the rest of this section, we reduce the above conjecture to a maximal rank problem. This maximal
rank problem for the particular case we need will be handled in the next section.

We want to apply the result of the previous section.
So we start from a flag $\HHHH  \subset \LLLL \subset \KKKK \subset \p^n$ where $\LLLL$ is a $t$-fat
$r$-plane with support $\HHHH$ and linear span $\KKKK$.
Our
first task is to identify the tangent space $T_\LLLL \FP$ at $\LLLL$ to the variety $\FP$ of fat planes. Recall that
the tangent space at $\LLLL$ to the full Hilbert scheme is
  $ H^0(\LLLL , N_\LLLL) $, where $N_\LLLL  := \mathcal{ H}om ( \mathcal{ I} _\LLLL, \OO_\LLLL)$ is the normal bundle.
Hence we look for a subspace
of that vector space.
We choose coordinates $x_i$ where $\LLLL$
is defined by the equations $x_0^t =x_1 = \cdots = x_{n-r-1} =0$ so
that we may identify $N_\LLLL$ as the direct sum $ \OO_\LLLL (1) ^{n-r-1}
\oplus \OO_\LLLL(t)$ and accordingly $H^0(N_\LLLL)$ as the direct sum $
H^0(\OO_\LLLL(1)) ^{n-r-1} \oplus H^0(\OO_\LLLL(t)).$

For the following lemma, we will introduce again a notation $\Ko$. The reader should be aware that, 
in the present section, this notation is introduced in such a way that $\Ko$ differs from $\Ho$ only in the special case
where $r$ is zero. For each integer $a$, we denote by $\Ko (\OO_\LLLL(a))$ the image of the restriction 
$ H^0 (\OO_{\p^n}(a)) \to H^0 (\OO_\LLLL(a) )$. We also extend this notation to sequences in the natural way:
by $\OO_\LLLL(\dd)$ we mean $\oplus_i \OO_\LLLL(d_i)$, and 
 $\Ko (\OO_\LLLL(\dd))$ stands for  $\oplus_i\Ko (\OO_\LLLL(d_i))$.
Finally we set $p := n-r-1$.

\bl \label {im1}

(i) The image of the natural morphism $ j : T\p^n \to  N_\LLLL$ from the tangent sheaf of $\pn$ to the normal sheaf $N_\LLLL$ is a subsheaf
$N'_\LLLL$ isomorphic to $ \OO_\LLLL(1)^{p} \oplus \OO_\HHHH (1)$ as an $\OO_\LLLL$-module.

(ii) More precisely, we may choose an isomorphism between $N_\LLLL$ and $ \OO_\LLLL(1)^{p} \oplus \OO_\LLLL (t)$ so that
the corresponding injection $ \OO_\LLLL(1)^{p} \oplus \OO_\HHHH (1)     \to \OO_\LLLL(1)^{p} \oplus \OO_\LLLL (t)$
decomposes as $\iota \oplus \mu$, where $\iota$ is the identity on the first summand and $\mu$
is the  multiplication by $x_0^{t-1}$ on the
second one.

(iii) The tangent space $T_\LLLL \FP$ at $\LLLL$ to the variety $\FP$ of fat planes is the image
of $H^0 (T\pn)$ (or $H^0 (\OO_{\pn} (1) ^{n+1})$) into $H^0(N'_\LLLL)$ (or $H^0(N_\LLLL)$). We write
$\Ko(N'_\LLLL)$ for this image.

(iv) 
Under the identification in (i), $\Ko(N'_\LLLL)$ appears as $
\Ko (\OO_\LLLL(1)) ^{p} \oplus H^0(\OO_\HHHH(1)).$
\el

\begin{proof} Let us start with the third statement. Since $\FP$ is the orbit in the full Hilbert scheme of 
$\LLLL$ under the projective
linear group, $T_\LLLL \FP$ has to be the image of the natural map $ H^0(j) : H^0(T\p^n) \to H^0(N_\LLLL).$

Now we turn to (i) and (ii).  Using our coordinates, our morphism $j$, viewed from
$ \OO_{\p^n}(1)^{n+1}$ to $ \OO_\LLLL(1)^{p} \oplus \OO_\LLLL (t)$ is given by the partial derivatives or our $n-r$ equations, which gives
essentially the announced matrix: just
note that, thanks to the characteristic zero assumption, the image
of the multiplication by the partial derivative $t x_0^{t-1}$ from $ \OO_\LLLL (1)$ to
$ \OO_\LLLL (t)$ is the same as the image of the multiplication by $x_0^{t-1}$, and this image is
 isomorphic to $ \OO_\HHHH (1)$.

Now we turn to (iv). We just note that by i), $H^0(N'_{\LLLL})$ is equal to 
$H^0(\OO_\LLLL(1))^{p}\oplus H^0(\OO_\HHHH(1))$. By definition,
$K^0(N'_\LLLL)$ is the image of $H^0 (\OO_{\pn} (1) ^{n+1})$ into $H^0(N'_\LLLL)$, which can now be identified as the space
$\Ko (\OO_\LLLL(1)) ^{p} \oplus H^0(\OO_\HHHH(1)).$

\end{proof}

Now we consider a flag $i : \LLLL \to Y$ of complete intersection subschemes in $\p^n$ where $Y$ is the general
complete intersection of type $\dd$ containing $\LLLL$. The tangent space at $Y$ to the corresponding
Hilbert scheme is $H^0(Y, N_Y)$ which can be computed as   $H^0(Y, \OO_Y(\dd) )$. It follows that
 $H^0(\LLLL, N_{Y|\LLLL})$ is isomorphic to $H^0(\LLLL, \OO_\LLLL (\dd))$.
Thus it is sound to write $\Ko(\LLLL, N_{Y|\LLLL})$ or simply  $\Ko( N_{Y|\LLLL})$ for the image
of $H^0 (N_Y)$ into
$H^0(\LLLL, N_{Y|\LLLL})$.

\bl \label{critabs}

(i) For $r>0$, $\Ko( N_{Y|\LLLL})$ is the whole of $H^0(\LLLL, N_{Y|\LLLL})$;

(ii) In any case, the dimension of  $\Ko( N_{Y|\LLLL})$ is $c$;

(iii) In any case, the natural map $ \Ko(N'_\LLLL) \to H^0 (\LLLL, N_{Y|\LLLL}) $ factors through $\Ko( N_{Y|\LLLL})$;

(iv) If the induced map
$ \Ko(N'_\LLLL) \to \Ko (N_{Y|\LLLL})$ is onto, then $\DD' \to \HHH'_2$ is smooth at $(\LLLL, Y)$.
\el
\begin{proof}
(i) We know that
$N_Y$ is
the direct sum $\OO_Y(\dd)$ and a standard cohomological argument
shows that all its sections come from the ambient projective
space. Hence what we have to prove is that any section of $\OO_\LLLL(\dd)$
comes from the ambient projective space, which follows from the
standard cohomological argument: the cohomology of line bundles on
projective spaces of dimension at least two is trivial (here we use $r
\ge 1$).

(ii) As  we have just seen, $\Ko( N_{Y|\LLLL})$ is the image of $H^0 (\OO_{\pn} (\dd))$ in $H^0 (\OO_\LLLL (\dd))$, so this is just
a count of monomials which we leave to the reader.

(iii) and (iv) Now we apply Prop.\ref {restr} : in our case, we have $\Ko (N_Y) = H^0(N_Y)$ and, according to Lemma~\ref{im1},
the assumption in Prop.\ref {restr} is precisely
the previous item. The statements (iii) and (iv) here are exactly the conclusions (ii) and (iv) there.

\end{proof}

We turn to the final result of the present section where, for sake of
clarity, we handle separately
the case $r=0$. 
We will write $\un, \deu$ and  $\ttt$ respectively for the
sequence $( 1, \cdots , 1), \,
(2 , \cdots, 2)$ and $(t, \cdots , t)$,
hence accordingly $\dd - \un, \, \dd - \deu, \, \dd - \ttt$ respectively for 
$(d_1 -1, \cdots, d_s -1), \, (d_1 -2, \cdots, d_s -2), \, (d_1 -t, \cdots, d_s -t).$

\bp \label{carz}

(i) For $r \ge 1$, we consider the generic morphism\\ $m : \OO_\LLLL
(1)^{p}\oplus \OO_\HHHH (1) \to \OO_\LLLL(\dd)$
of coherent $\OO_\LLLL$-modules. If 
$$
H^0(m) : H^0 (\OO_\LLLL (1)^{p}) \oplus H^0 (\OO_\HHHH (1)) \to
H^0 ( \OO_\LLLL(\dd))
$$
is onto then $\DD' \to \HHH'_2$ is onto too.

(ii) For $r = 0$, we denote by
$\Ko (\OO_\LLLL (1)^{p}\oplus \OO_\HHHH (1), \OO_\LLLL(\dd))$
the image of
the natural map 
$$
 \Ko( \OO_{\LLLL}(\dd - \un) ) \oplus
\Ko( \OO_{\LLLL}(\dd - \ttt ) ) \rar Hom ( \OO_\LLLL (1)^{p}\oplus \OO_\HHHH
(1),  \OO_\LLLL(\dd) ).
$$
We consider the generic morphism $m$ in $\Ko (\OO_\LLLL (1)^{p}\oplus
\OO_\HHHH (1), \OO_\LLLL(\dd))$.
If the image by $H^0(m)$ of
$\Ko (\OO_\LLLL (1)^{p})\oplus H^0(\OO_\HHHH (1))$ into $H^0( \OO_\LLLL(\dd))$
is $\Ko ( \OO_\LLLL(\dd)),$ then $\DD' \to \HHH'_2$ is onto.
\ep

\begin{proof}
(i) Since our morphism $\DD' \to \HHH'_2$ is projective and the codomain is irreducible, it is sufficient 
to prove that it is dominant.
We apply Lemma \ref{critabs} (iv), hence we have to prove that the map
$ m_Y : H^0(N'_\LLLL) \to H^0 (\LLLL, N_{Y|\LLLL})$ is onto. This map depends upon
our complete
intersection $Y$. We express it in terms of
the system of equations $b := (b_1, \cdots , b_s) \in H^ 0 ({\mathcal I}_\LLLL (\dd))$ of
$Y$, rather than in terms of $Y$
itself.
This allows us to describe the associated morphism $m_Y : N'_\LLLL \to
N_{Y|\LLLL}$
or, via the identifications of Lemma~\ref{im1},
$m_b : \OO_\LLLL (1)^{p}\oplus \OO_\HHHH (1) \to
\OO_\LLLL(\dd)$ as follows:

- for the first factor, the $j$-th component ($1 \le j \le p$), from
$ \OO_\LLLL (1)$ to
$\OO_\LLLL(\dd)$, is the derivative of $b$ with respect to $x_j$;

-for the second factor, from $ \OO_H (1)$ to $\OO_\LLLL(\dd)$, we have the
derivative of $b$ with respect to $x_0$ (note that indeed this
derivative factors through $ \OO_H (1)$).

What we have to prove is that, for $b$ sufficiently general, $\Ho(m_b)$ is onto.
For this, thanks to our surjectivity assumption, it is enough to prove that
 $b \mapsto m_b $ is dominant (or onto).

We prove that $b \mapsto m_b $  is onto. For this we take $m :=
(m_1, \cdots , m_{p}, m_0) : \OO_\LLLL (1)^{p}\oplus \OO_\HHHH (1)
\to \OO_\LLLL(\dd)$ and search for $b$ with $m = m_b$.  By the standard
cohomological argument, we may lift $m_1, \cdots , m_{p}$ and
consider we are given $(m_1, \cdots , m_{p}) : \OO_{\pn} (1)^{p}
\to \OO_{\pn}(\dd)$. Now for $m_0$, we see it as a section of $Hom
(\OO_\HHHH , \OO_\LLLL( \dd - \un) )$ hence as a section of $ \OO_\LLLL( \dd -
\un)$ annihilated by $x_0$, thus of the form $tx_0^{t-1}f$ with $f$ a
section of $\OO_\LLLL( \dd - \ttt)$, using the characteristic zero
assumption.  As above, we may lift $f$ as a section, still denoted $f$
of $\OO_{\pn} ( \dd - \ttt)$. At this point we may set $b:= x_1 m_1 +
\cdots + x_{p} m_{p} + x_0^t f$ and check that it has the
desired property.

(ii) The proof is almost the same: we apply Lemma \ref{critabs} (iv). This
time, we have to prove
that the map $ m_Y : H^0(N'_\LLLL) \to H^0 (\LLLL, N_{Y|\LLLL})$ sends
$ \Ko(N'_\LLLL)$ onto $\Ko (\LLLL, N_{Y|\LLLL})$. As above we introduce
a system of equations $b := (b_1, \cdots , b_s) \in H^ 0 ({\mathcal I}_\LLLL (\dd))$ of
$Y$.
Via the identifications of Lemma~\ref{im1}, we are concerned, for $b$
sufficiently general, by the sheaf morphism
$m_b : \OO_\LLLL (1)^{p}\oplus \OO_\HHHH (1) \to
\OO_\LLLL(\dd)$ defined by the same formulas as in the previous case.

Since our identifications send  $ \Ko(N'_\LLLL)$ to $\Ko (\OO_\LLLL
(1))^{p}\oplus  H^0( \OO_\HHHH (1))$ and $\Ko (\LLLL, N_{Y|\LLLL})$ to $\Ko
(\OO_\LLLL(\dd))$,
it is enough to prove that the
image  of $b \mapsto m_b $ is $\Ko (\OO_\LLLL (1)^{p}\oplus \OO_\HHHH
(1), \OO_\LLLL(\dd))$.
For this we take $m :=
(m_1, \cdots , m_{p}, m_0)$ in $\Ko (\OO_\LLLL (1)^{p}\oplus \OO_\HHHH
(1), \OO_\LLLL(\dd))$.
This means that $m_1, \cdots , m_{p}$ come from  sections still denoted
$m_1, \cdots , m_{p}$
in $\Ho(\pn, \OO (\dd - \un))$, while $m_0$ is of the form $x_0^{t-1}f$, or
better of the form
$tx_0^{t-1}f$,
with $f$ a
section of $\Ho (\pn, \OO( \dd - \ttt))$.
And we search for $b$ with $m = m_b$. Again we may set $b:= x_1 m_1 +
\cdots + x_{p} m_{p} + x_0^t f$ and check that it has the
desired property.

\end{proof}


\section{Nonemptiness}
\label{maxrank}

In this section, we prove our conjecture \ref{cdimfat} in the case we need.
We restrict to the very special case where $t$ is the greatest number in our sequence $\dd$,
and we assume furthermore that $t$ is at least $3$, and that it occurs only once in $\dd$.
We will prove:

\bp \label{nonempty}
Under the above restrictions, when
$\rho$
is nonnegative, for any complete intersection $Y$ of type $\dd$ in $\p^n$, 
the Hilbert scheme of
$t$-fat $r$-planes in $Y$ is nonempty.
\ep

We keep the notations of the previous section. Furthermore, we
denote by $h^0(u, e)$ the number of monomials of degree
$e$ in $u$ variables,
and accordingly, for any sequence $\ee := (e_1, \cdots , e_s)$ of integers, we set
$ h^0(u, \ee) := h^0(u, e_1)+ \cdots + h^0(u, e_s).$

Thanks to Proposition \ref{carz}, it is enough to prove a
maximal rank statement, which depends on whether $r$ is zero or not.
Namely, we have to prove the following two lemmas.

\bl \label {lmaxrank} (the case $r \ge 1$)
For $p$ satisfying $(r+2)p +r+1 \ge h^0(r+2, \dd) -1$, and
for the general morphism $m : \OO_\LLLL (1)^{p}\oplus \OO_\HHHH (1) \rar \OO_\LLLL
(\dd)$, $H^0(m)$ is onto.
\el

\bl \label {lmaxrankz} 
(the case $r = 0$) We suppose $2p +1 \ge h^0(2, \dd) -1$.
Then,
for the general morphism $m$ in $\Ko (\OO_\LLLL (1)^{p}\oplus
\OO_\HHHH (1), \OO_\LLLL(\dd))$, the image by $H^0(m)$ of
$\Ko (\OO_\LLLL (1)^{p})\oplus H^0(\OO_\HHHH (1))$ into $H^0( \OO_\LLLL(\dd))$
is $\Ko ( \OO_\LLLL(\dd)).$ 

\el

Recall that 
$\Ko (\OO_\LLLL (1)^{p}\oplus \OO_\HHHH (1), \OO_\LLLL(\dd))$ denotes
the image of
the natural map from $ \Ko( \OO_{\LLLL}(\dd - \un) ) \oplus
\Ko( \OO_{\LLLL}(\dd - \ttt ) )$ to   $Hom ( \OO_\LLLL (1)^{p}\oplus \OO_\HHHH
(1),  \OO_\LLLL(\dd) ).$

The differences between our two lemmas can be erased by switching to the point 
of view of graded modules.
So, just for the present section, we radically change the meaning of our notations: 
from now on, $\OO_\LLLL$ denotes the graded ring $k[x_0,  \cdots , x_{r+1}]/(x_0^t)$ and
$\OO_\HHHH$ denotes the quotient graded module $k[x_1,  \cdots , x_{r+1}]$.
For a graded $\OO_\LLLL$-module $G$ with graduation $\gamma$, we write $G(a)$ for
the module $G$ equipped with the graduation $\gamma_a := \gamma-a$. 
For a graded-module $G$, by $\Ho(G)$ we mean the degree $0$ component of $G$, while for
a morphism $m$ of graded modules, by $\Ho (m)$ we mean the restriction of $m$ to
the degree $0$ components.
With these conventions, our two lemmas rephrase as the single following one:

\bl \label {lmaxrankg} 
For $p$ satisfying $(r+2)p +r+1 \ge h^0(r+2, \dd) -1$, and
for the general morphism $m : \OO_\LLLL (1)^{p}\oplus \OO_\HHHH (1) \rar \OO_\LLLL
(\dd)$, $H^0(m)$ is onto.
\el

Proof.
Here we use a method which can be tracked back to \cite {EH, EHM}, where similar results were 
obtained in a different context.
We denote by $M$ the vector space $Hom (\OO_\LLLL (1)^{p}\oplus \OO_\HHHH
(1), \OO_\LLLL (\dd))$, by $S$ the
space of nontrivial linear forms on $\Ho (\OO_\LLLL(\dd))$, and by $\ZZ$ the
``incidence'' subscheme in $M \times S$ consisting of pairs $(m, \ell)$
for which $\ell \circ \Ho(m)$ vanishes. We denote by $\ell$ the second projection:
$\ell : \ZZ \to S$. What we want to prove is that the first projection
$\ZZ \to M$ is not dominant. This will follow if we prove
the inequality $dim \ZZ \le dim M$, since the fibers of our projection are unions of lines. 
We proceed by contradiction and
suppose that the projection
$\ZZ \to M$ is dominant. To each $\lla \in S$ we attach the bilinear form
$\lla^*$ on $\Ho ( \OO_\LLLL (\dd -\un)) \times \Ho ( \OO_\LLLL (1)) $ defined by
$\lla^*( f, v)= \lla (vf)$.
By semi-continuity, we have an open
subset $\ZZ_u \subset \ZZ$ which still dominates $M$, and where the rank $u$
of $\ell^*$ is constant.

In the first factor $\Ho ( \OO_\LLLL (\dd -\un)) $ of our product, we have a
distinguished line:
the line $D$ generated (in the summand $\OO_\LLLL (t)$) by
$x_0 ^{t-1}$. Our first observation is the following:

\bsl \label {vanish}

For our general point $z \in \ZZ_u$, $\lls (z)$ vanishes on $D \times \Ho (
\OO_\LLLL (1)) $.
\esl
\begin{proof}
 In $ \OO_\LLLL (\dd ) $ we have the summand $\OO_\LLLL (t)$.
And therein, we have the graded submodule $x_0^{t-1}\OO_\LLLL (1) $ consisting of multiples of
$x_0^{t-1}$. This submodule may be better denoted by $x_0^{t-1}\OO_\HHHH  (1)  $ 
since the multiplication by $x_0^{t-1}$, which sends $\OO_\LLLL  (1)   $ into
$\OO_\LLLL (t)  $, factors through $\OO_\HHHH  (1) $.  This submodule $x_0^{t-1}\OO_\HHHH  (1)  $ 
is easily identified as
the submodule of $ \OO_\LLLL (\dd ) $ which is annihilated by $x_0$.
Hence, any morphism $m \in M$ has to send the second summand $\OO_\HHHH
(1) $ of its domain, which is
annihilated by $x_0$,
into the summand $\OO_\LLLL (t)$ of its codomain, and more precisely into
the submodule $x_0^{t-1} \OO_\HHHH (1)$ mentioned above.
Also a sufficiently general morphism $m \in M$ sends $\OO_\HHHH (1) $
isomorphically onto that submodule.
Accordingly, $\Ho(m)$ sends $\Ho(\OO_\HHHH (1)) $ isomorphically onto
$\Ho(x_0^{t-1} \OO_\HHHH (1))$.
So, for our general
$z \in \ZZ_u$, $\ell (z)$ has to vanish on $\Ho(x_0^{t-1} \OO_\HHHH (1))$
which implies that $\lls (z)$ vanishes on
$D \times \Ho ( \OO_\LLLL (1)) $.

\end{proof}

Our next observation
stresses the role of $u$, which is to control
the dimension of the fiber of $\ZZ \to S$.
We denote by $S_u$ the projection of $\ZZ_u$ in $S$.

\bsl \label {dimfib}
The codimension of the fiber of $\ZZ_u$ over a point $\lla \in S_{u}$ is
$pu$.
\esl

Proof:
Let $m := (m_1, \cdots , m_{p}; m_0) $ be a point in $M$, where $(m_1,
\cdots , m_{p})$ are in
$\Ho( \OO_\LLLL (\dd -\un))$ while $m_0$ is in $H^0( \OO_\LLLL (\dd -\ttt))$.

Thanks to the previous lemma,
we see that $\lla \circ \Ho(m)$
vanishes if and only if $\lla^* (m_1) = \cdots =\lla^* (m_{p})=0.$ 
Each one among these $p$ equations imposes $u$ independent conditions on
$m$, since the rank of $\lla ^*$ is $u$. Since these
equations concern different components of $m$, their ranks add up to the
rank of $m \mapsto \lla \circ \Ho(m)$
which turns out to be $pu$.

\qed

Our next task consists in estimating the dimension of $S_u$.

\bsl \label {dimSuu}
The dimension of $S_u$ is at most $ \ho(u, \dd)+ (r+1 -u)u .$
\esl
\begin{proof}
For this we have to single out the line $E$ generated by $x_0$ in
$\Ho ( \OO_\LLLL (1)) $ and to distinguish two cases
according to whether, for our general $z \in \ZZ_u$, $\lls(z)$ vanishes or
not on
$\Ho ( \OO_\LLLL (\dd -\un)) \times E $.

(i) We start with the (slightly simpler) case where $\lls(z)$ does not
vanish on\\
$\Ho ( \OO_\LLLL (\dd -\un)) \times E $.

In order to bound  the dimension of $S_{u}$  at a point $\lla_0$, we will
define, in a neighborhood $U \subset S_u$
of $\lla_0$,
two algebraic
maps $f : U \to \AAA^b$ and $g : U \to \AAA^c$ so that $(f, g)$ is injective.
This will bound the
dimension of $S_u$ by $b+c$.
To this effect, we reorder our basis $C := (x_0, \cdots , x_{r+1})$ of $ \Ho(
\OO_{L'} (1)) $ (where $x_0$ remains an equation
of $\HHHH$) so that, in this basis, the first $u$ rows of the
matrix of $\lla ^*_0$ are linearly independent.
This property will hold in a neighborhood of $\lla _0$ which we take as $U$.
We write $C'$ for the sub-basis $(x_0, \cdots , x_{u-1})$ and $C''$ for the
rest of the basis so that we have
$C = C' \amalg C''$. Next, in $ \OO_\LLLL (\dd)$, we have the basis consisting
of monomials in each summand, which we call $\dd$-monomials. Similarly, we have
the basis of $(\dd - \un)$-monomials in  $ \OO_\LLLL (\dd- \un)$.

Associated with these bases,
we have the matrix $N_\lla$ of $\lla^*$, which is an algebraic
function of $\lla$.
Now for each element in $C''$,
we have the $u$ coordinates of the corresponding row (in $N_\lla$) as a
combination of the rows in $C'$.
This defines
$b:= (r+2 -u)u$ functions on $U$ which altogether yield our map $f$.

Now for $g (\lla)$ we take the restriction $\lla'$ of
$\lla$ to the subspace
generated by the following set $T'$ of $\dd$-monomials : at first, take the
set $T$ of $\dd$-monomials
depending (at most) on variables in $C'$, then delete those, in the 
summand $\OO_\LLLL (t)$, which
are divisible by $x_0^{t-1}$. This deletion corresponds to the fact observed
above that $\lla$ vanishes there.

What we have to check is that $\lla$ is
determined by $\lla'$ and $f(\lla)$.
For this, we claim that for each integer $q$ with $0 \le q \le t$ 
the values of $\lla$ on the set $T_q$ of those
$\dd$-monomials which are
of degree $q$ with respect to variables in $C''$ are linear combinations (where
coefficients are polynomials in $f(\lla)$)
of its values on $T'$. We prove the claim by induction on $q$, the case $q := 0$ following from the
vanishing mentioned above.
For the general case we consider a $\dd$-monomial $m := m'x_i$ where $m'$ is
a $(\dd - \un)$-monomial
and $x_i$ is in $C''$. In the column corresponding to $m'$ in $N_\lla$, the
first $u$ entries are values of $\lla$
on elements of $T_{q-1}$, while the entry in the row corresponding to $x_i$
is $\lla (m)$, which gives us the
desired linear relation.

It remains to check that the number of elements in $T'$ is
$ \ko (u, \dd)-u$. Indeed, $1$ is subtracted from $ \ko (u, \dd)$ 
because, although $\dd$ contains
$t$, $x_0 ^t$ is not a
$\dd$-monomial,
and $u-1$ is subtracted due to the difference between $T$ and $T'$. Thus the codomain of
our map
$g$ is $\AAA^c$ with
$c := \ho (u, \dd)-u.$

(ii) Now we treat the similar case where $\lls(z)$ vanishes on
$\Ho ( \OO_\LLLL (\dd -\un)) \times E $. The method is the same so we just
highlight the changes.
Thanks to the vanishing assumption, $\lla^*$ is now determined by the bilinear
form $\lla'^*$ induced on
$\Ho ( \OO_\LLLL (\dd -\un)) \times \Ho ( \OO_\HHHH (1)) $. Our basis $C$ now
has the form
$(x_1, \cdots , x_{r+1})$, and the subbasis
$C'$ is $(x_1, \cdots , x_{u})$  Accordingly, the number $b$ is now equal to
$(r+1 -u)u$. On the other hand,
here, there is no deletion, $T'$ is equal to $T$ and its number of elements
is $ \ho (u, \dd),$ which yields
the desired formula.
\end{proof}

In order to complete the proof of \ref{lmaxrank}, it remains to check that
the estimates obtained so far
make the dimension of $\ZZ_u$ smaller than that of $M$, namely that the
codimension (in $M$) obtained
for the fiber of $\ZZ_u \to S_u$ is bigger than the dimension of $S_u$. This
reads:

\bsl \label {lastineg}

For $t \ge 3$, $p$ satisfying $(r+2)p +r+1 \ge \ho(r+2, \dd) -1$, and $1 \le
u \le r+2$ we have
$ \ho(u, \dd)+ (r+1 -u)u \le pu$.
\esl
\begin{proof}
We argue by convexity (with respect to $u$) and start by checking the
extreme cases:

(i)
For $u := r+2$, the desired conclusion is just the assumption.

(ii)
For $u:=1$, we contrapose and prove that $p \le \ho (1, \dd) +r -1$ implies\\
$(r+2)p +r+3 \le \ko(r+2, \dd). $
Taking the critical value $s+r-1$ for $p$
we have to prove
$s(r+2)+(r+1)^2  \le \ho(r+2, \dd). $

We split this inequality summand by summand, in other words we claim

a) $r+2 +(r+1)^2  \le \ho(r+2, t)$ (for the occurence of $t$ in $\dd$) and

b) $r+2  \le \ko(r+2, \delta) $ (for each other integer, $\delta \ge 2$, in $\dd$).

For a) it is sufficient to check the first case $t :=3$. In this case, we
have to prove
$6(r+1)^2  \le (r+2)[(r+3)(r+4) -6]$  or, dividing by $r+1$,
$6r+6 \le (r+2)(r+6)$, or $0 \le r^2+2r+6$, which is evident. While b) is
clear since for each variable $x_i$,
we have the monomial $x_i^\delta$.

It remains to check that the function $f := u \mapsto \ho(u, \dd)+ (r+1 -u)u
$ is convex on our interval
$[1, r+2]$. For this, we compute the discrete derivatives $f' := u \mapsto
f(u+1) -f(u)$ and $f''$. We find
$f'(u) = \ho(u+1, \dd - \un) +  r + 1 -2u -1$ and $f''(u) =  \ho(u+2, \dd - \deu)
-2$.
We see that
this second derivative is nonnegative for $u \ge 1$, yielding the desired
convexity.
\end{proof}


\section{Spannedness}\label{spannedness}

This section is devoted to the proof of the desired covering statement :

\bp
\label{strongcov}
Let $r \ge 0$, $\, \, 1 \le s \le n-r-1$ and $2 \le d_1 \le ... \le d_{s-1} < d_s$ be
integers. We set
$\rho :=
(r+2)(n-r)-\Sigma_{i=1}^s {d_i + r +1 \choose r+1}$ and assume 
the (necessary) inequality

$$
\rho +r \ge n-s.
$$
If $Y \subset Y'$ is any pair of type $(d_1, ..., d_s)$ in $\pn$, then $Y$
is covered
by strong $r$-planes.
\ep

Thanks to Proposition \ref{fatstrong},  this statement is an immediate consequence of the following one.

\bp \label {fatcov}
Under the same assumptions,  $Y$ is covered
by $d_s$-fat $r$-planes.
\ep

We pose $t := d_s$.

From \S \ref{fatplanes}, we have the restricted flag-Hilbert scheme $\DD' \subset
\HHH'_1 \times \HHH_2$.
Over the first factor $\HHH'_1$, we have the universal $t$-fat $r$-plane, say $\LL \subset \HHH'_1 \times \pn$.
Over the second factor $\HHH_2$, we have the universal complete
intersection of type $\dd$, say $\UU \subset \HHH_2 \times \pn$,
and over $\DD'$, we have the universal flag, say $\hLL \subset \hUU
\subset \DD' \times \pn$.
We have a natural projection $e : \hUU \to \UU$, and what we have to prove is that the
restriction $e' : \hLL \to \UU$ is onto.

Since
$e'$ is a $\HHH_2$-morphism among varieties which are
projective over $\HHH_2$, its image is also projective over $\HHH_2$. So it
is sufficient to prove that $e'$
is dominant, and, for that, to find one point in $\hLL$ where the fiber of
$e'$
has the expected dimension $\rho +r -n +s$, and not more.

So we compute the fiber of $e'$ at a point $(L, Y, p) \in \hLL$, where $L$
is a $t$-fat $r$-plane contained in the complete intersection $Y$ and $p$ is a point on $L$.
This splits into two cases according to whether $r$ is zero or not.

(i) The case $r=0$. This case is known since \cite{Ro}. Hence we just give the idea of the proof, which is similar but simpler 
than the other case. The variety of $t$-fat points at $p$ contained in $Y$
is identified with a subvariety in
the projectivized tangent space of $Y$ at $p$ with equations depending on the equations of $Y$.
The number of these equations is easily checked
to be 
$\Sigma_i (d_i -1) -1$: $d_i -1$ is the number of degrees between $2$ and $d_i$,
and $1$ is subtracted for the degree $t$. Thanks to our assumption on $\rho$, this is
at most $n-s -1$ which is the dimension of this projective space. Hence this variety is nonempty.

(ii) The case $r \ge 1$.
We consider the subscheme
$\cV_{Y,p}$
in the projectivized tangent space $PT_{Y,p}$ of $Y$ at the point $p$, which
parametrizes lines through $p$ contained in
$Y$. It is defined by the homogeneous components of the Taylor expansions at
$p$ of the
equations of $Y$. For this reason,  we set\\
$\dd':= (1, 2, \cdots, d_1, \cdots , 1, 2, \cdots , d_s)$ 
and $\dd'':= (2, \cdots, d_1, \cdots , 2, \cdots , d_s)$.

We immediately observe that there is a natural isomorphism between the
vector space of tuples of equations of
type $\dd$ (in $n+1$ variables) vanishing at $p$ and tuples of equations of
type $\dd'$ (in $n$ variables).
This yields:

For $(Y,p)$ sufficiently general in $\UU$, $\cV_{Y,p}$ is a sufficiently
general
complete intersection of type $\dd''$ in $PT_{Y,p}$.

Next we have the following

\bl \label{fibre}
For any $(Y, p) \in \UU$, the fiber of $e'$ over $(Y,p)$ is isomorphic with
the Hilbert scheme of $t$-fat
$(r-1)$-planes in $\cV_{Y,p}$.
\el

First we check how we may complete the proof of Proposition \ref{fatcov} using this lemma. Thanks to 
Proposition \ref{nonempty}, we just have to check
that the expected dimension $\rho'$ of the
Hilbert scheme of $t$-fat $(r-1)$-planes in the generic complete
intersection of type $\dd'$ in $\p^{n-1}$ or equivalently of type $\dd''$ in
$\p^{n-s-1}$ is nonnegative. We have
\begin{eqnarray*}
\rho' & = & (r+1)(n-r)-1 -\Sigma_{i=1}^s\Sigma_{j=1}^{d_i} {j + r \choose r}\\
      & = &(r+1)(n-r)-1-\Sigma_{i=1}^s ({d_i + r +1 \choose r+1} -1) \\
      & = &(r+2)(n-r) -1 -n +r +s -\Sigma_{i=1}^s {d_i + r +1 \choose r+1} \\
      & = & \rho +r -n +s.
\end{eqnarray*}      
This is nonnegative by assumption.

Now we prove lemma \ref{fibre}.
First of all, we have a natural isomorphism $g$ between the Hilbert scheme
of $t$-fat $r$-planes in $\p^n$
passing through $p$ and the Hilbert scheme of $t$-fat $(r-1)$-planes in the
projectivized tangent space
$PT\p^n_{p}$: if we identify this projectivized tangent space with a
hyperplane $K \subset \p^n$ not passing
through $p$, $g(L)$ is the scheme-theoretic intersection of $L$ with $K$.

Now we prove that $g$ induces a bijection from the Hilbert scheme $H'$  of
$t$-fat $r$-planes in $Y$ passing
through $p$ to the Hilbert scheme $H''$ of $t$-fat $(r-1)$-planes in
$W_{Y,p}$. Let $M$ be a $t$-fat
$(r-1)$-plane in
$W_{Y,p}$ and $M'\subset \p^n$ be the unique $t$-fat $r$-plane through $p$
corresponding to $M$ (hence $g(M')=M$).
Let $(f_i := \Sigma_{j \le d_i} f_{ij}) _{i \le s}$ be a system of equations
of $Y$. Note that the $f_i$'s
vanish at $p$. Our claim is that the $f_i$'s all vanish on $M'$ if and only
if the $f_{ij}$
all vanish on $M=g(M')$. This follows readily from the particular case of a
single equation $f$, with $n=r+1$,
which we state explicitly:

Let $f(x_0, x'):= x_0f_1(x') + \cdots + x_0^\delta f_\delta (x')$ be a homogeneous
polynomial in $n+1$ variables
$x_0, \cdots, x_n$ where $x'$ stands for $(x_1, \cdots, x_n)$. Note that $f$
vanishes at $p:=(1, 0, \cdots,0)$.
We denote by $K$ the hyperplane defined by $x_0= 0$, by
$M$ the $t$-fat $(n-1)$-plane defined by $x_1^t =0$ and by $N$ the $t$-fat
$(n-2)$-plane defined
in $K$ by the same equation. Since $f$ is a multiple of $x_1^t$ if and only if all the $f_i$'s are,
we have that
$f$ vanishes on $M$ if and only if the
$f_i$'s vanish on $N$.

\qed

\end{document}